\documentclass[preprint,12pt]{elsarticle}

\usepackage{epsfig}

\usepackage{amssymb}
\usepackage{amsthm}

\journal{Journal of Geometry and Physics}

\newcommand{\eh}{\hspace{.05in}}
\newcommand{\ih}{\'\i}
\newcommand{\R}{\mathbb{R}} 
\newcommand{\A}{\mathcal{A}}
\newcommand{\C}{\mathbb{C}}
\newcommand{\E}{\mathbb{E}}
\newcommand{\CC}{\mathcal{C}}
\newcommand{\D}{\mathcal{D}}
\newcommand{\m}{\frac{1}{2}}   
\newcommand{\q}{{\frac{i\pi}{4}}}
  
\newcommand{\ds}{\displaystyle}  
  
\newcommand{\BE}{\begin{equation}}  
\newcommand{\EE}{\end{equation}}  
\newcommand{\Lim}[1]{\lower5pt\hbox{${{\ds\lim}\atop #1}$}}

\begin{document}

\begin{frontmatter}
\author{M.F. da Silva\corref{cor1}\fnref{label2}}
\ead{marcio.silva@ufabc.edu.br}
\ead[url]{http://sites.google.com/site/marcfab}

\author{G.A. Lobos\fnref{label3}}
\ead{lobos@dm.ufscar.br}
\ead[url]{http://www2.dm.ufscar.br/$\sim$lobos}

\author{V. Ramos Batista\fnref{label2}}
\ead{valerio.batista@ufabc.edu.br}
\ead[url]{http://valerioramosbatista.googlepages.com}

\fntext[label2]{Centre for Mathematics, Computer Science and Cognition, ABC Federal University, r. Catequese 242, 09090-400 \ St.Andr\'e-SP, Brazil}
\fntext[label3]{Mathematics Department, Federal University of S\~ao Carlos, rod. Washington Lu\ih s km 235, 13565-905 \ S\~ao Carlos-SP, Brazil}

\cortext[cor1]{Corresponding author.}

\title{Minimal surfaces with only horizontal symmetries}

\begin{abstract}
We find the first examples of triply periodic minimal surfaces of which the intrinsic symmetries are all of horizontal type.
\end{abstract}

\end{frontmatter}

\section{Introduction}
During the {\it Clay Mathematics Institute 2001 Summer School on the Global Theory of Minimal Surfaces}, M. Weber introduced the following definitions in his first lecture entitled {\it Embedded minimal surfaces of finite topology}:

``A {\it horizontal symmetry} is a reflection at a vertical plane or a rotation about a horizontal line. A {\it vertical symmetry} is a reflection at a horizontal plane or a rotation about a vertical line.'' 

With these definitions, he proved that such symmetries induce symmetries in the {\it cone metrics} $dh$, $Gdh$ and $dh/G$ for a {\it Weierstra\ss \ pair} $(G,dh)$ of a minimal surface (see \cite{[W1]} and \cite{[W2]} for details).  

By classifying the symmetries this way, we sort out the space groups that might admit one, both or none of them. Since minimal surfaces may model some natural structures, like crystals and co-polymers, an example within a given symmetry group might fit an already existing compound, or even hint at non-existing ones. However, several symmetry groups are not yet represented by any minimal surface (see \cite{[H]} and \cite{[LM]} for details and comments). 
   
Restricted to intrinsic symmetries, outside the triply periodic class it is easy to give examples of complete embedded minimal surfaces in $\R^3$ of which the symmetries are either only horizontal or vertical. For instance, the Costa surface (see \cite{[C]} and \cite{[K]}) has only horizontal intrinsic symmetries. The doubly periodic examples found by W. Meeks and H. Rosenberg in \cite{[MR]} have only vertical symmetries (see also \cite{[PRT]} for nice pictures).     

Curiously, in the widest class of minimal surfaces, namely the triply periodic, each example known before presents either both or none of such symmetries, according to our analysis. In this present work, we show examples which are probably the first triply periodic minimal surfaces with only horizontal symmetries. Regarding examples with only vertical symmetries, we believe they have not yet been found. 

The examples presented herein are inspired in the surfaces C$_2$ and L$_{2,4}$ from \cite{[V1]} and \cite{[V2]}. Any of those is generated by a fundamental piece, which is a surface with boundary in $\R^3$ with two catenoidal ends. The fundamental piece resembles the Costa surface with its planar end replaced by either symmetry curves or line segments. By suppressing the catenoidal ends, if we pile up several copies of the fundamental piece, we get the pictures in Figure 1 and Figure 2(b). They are also named C$_2$ and L$_{2,4}$. 

\input epsf  
\begin{figure} [ht] 
\centerline{ 
\epsfxsize 6cm  
\epsfbox{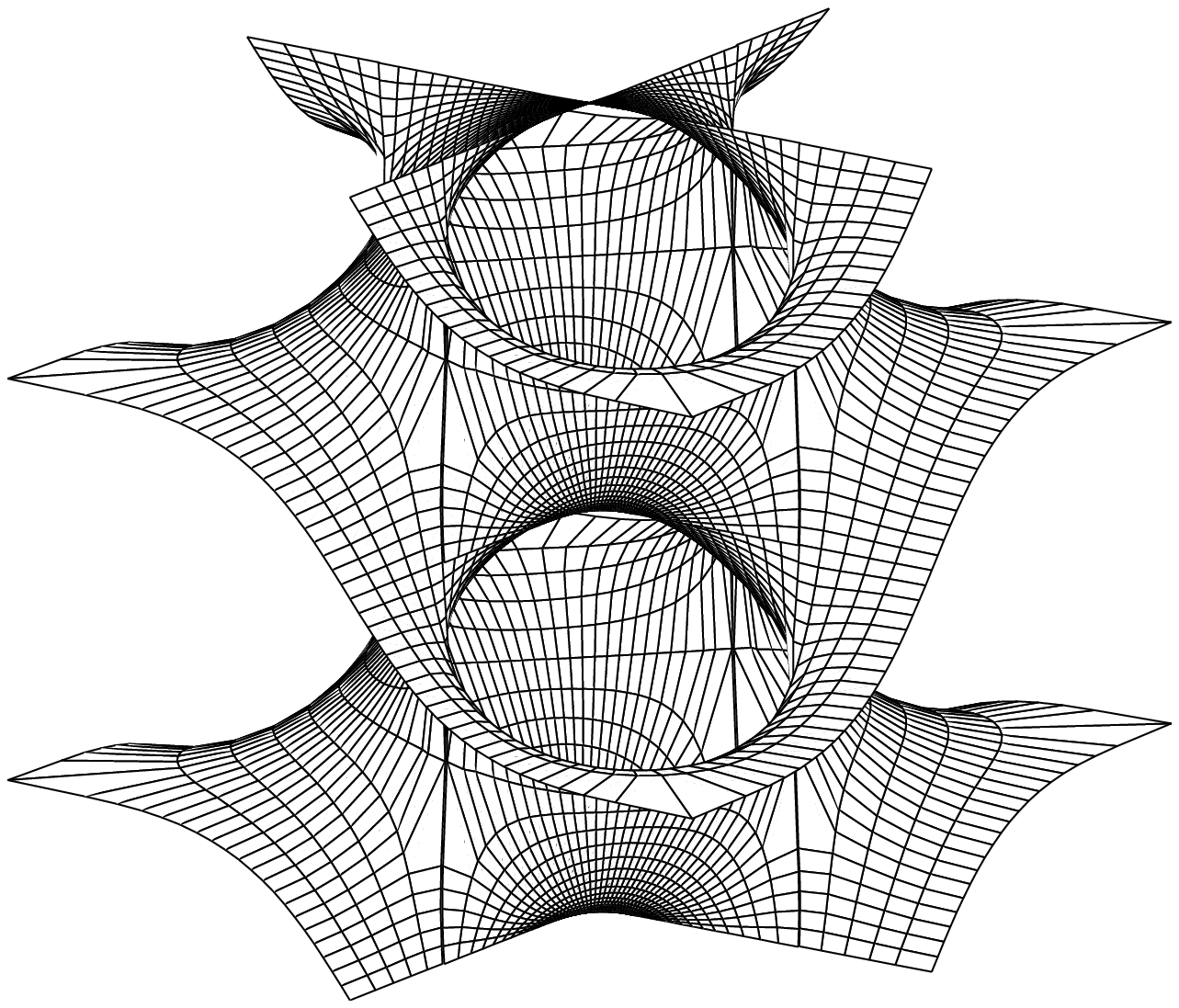}
\epsfxsize 7.5cm  
\epsfbox{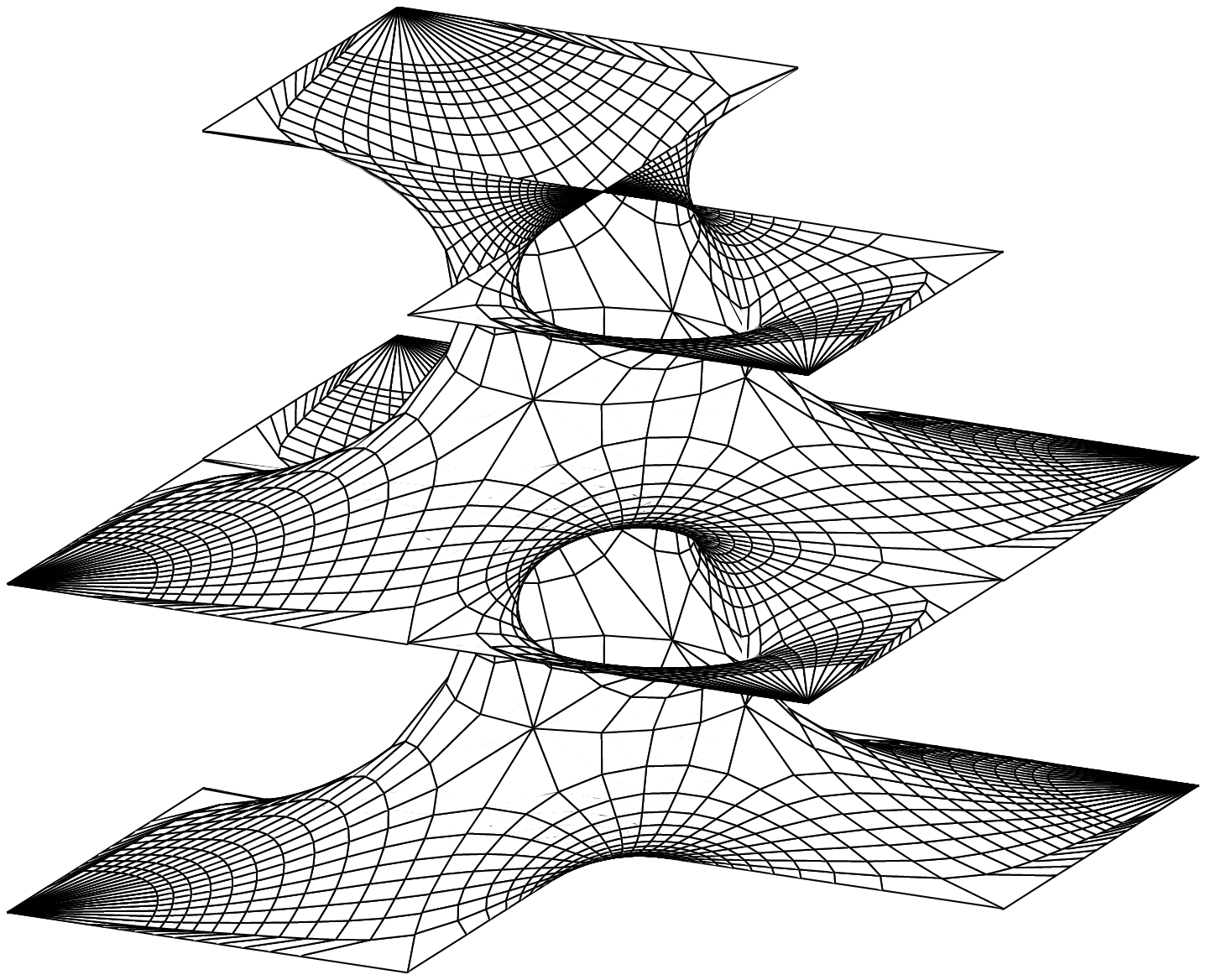}}
\centerline{(a)\hspace{7cm}(b)}  
\caption{(a) The surfaces C$_2$; (b) the surfaces L$_2$.} 
\end{figure} 
 
The reader will notice that the surfaces C$_4$, also described in \cite{[V1]} and \cite{[V2]}, were not mentioned beforehand. This is because, for them, the ``piling up'' procedure naturally forces extra symmetries to exist, and one goes back to a famous surface from H. Schwarz (see Figure 2(a)). Notice, for instance, the vertical straight line that comes out in the surface.

\input epsf
\begin{figure} [ht] 
\centerline{
\hspace{-3cm}
\epsfxsize 8cm  
\epsfbox{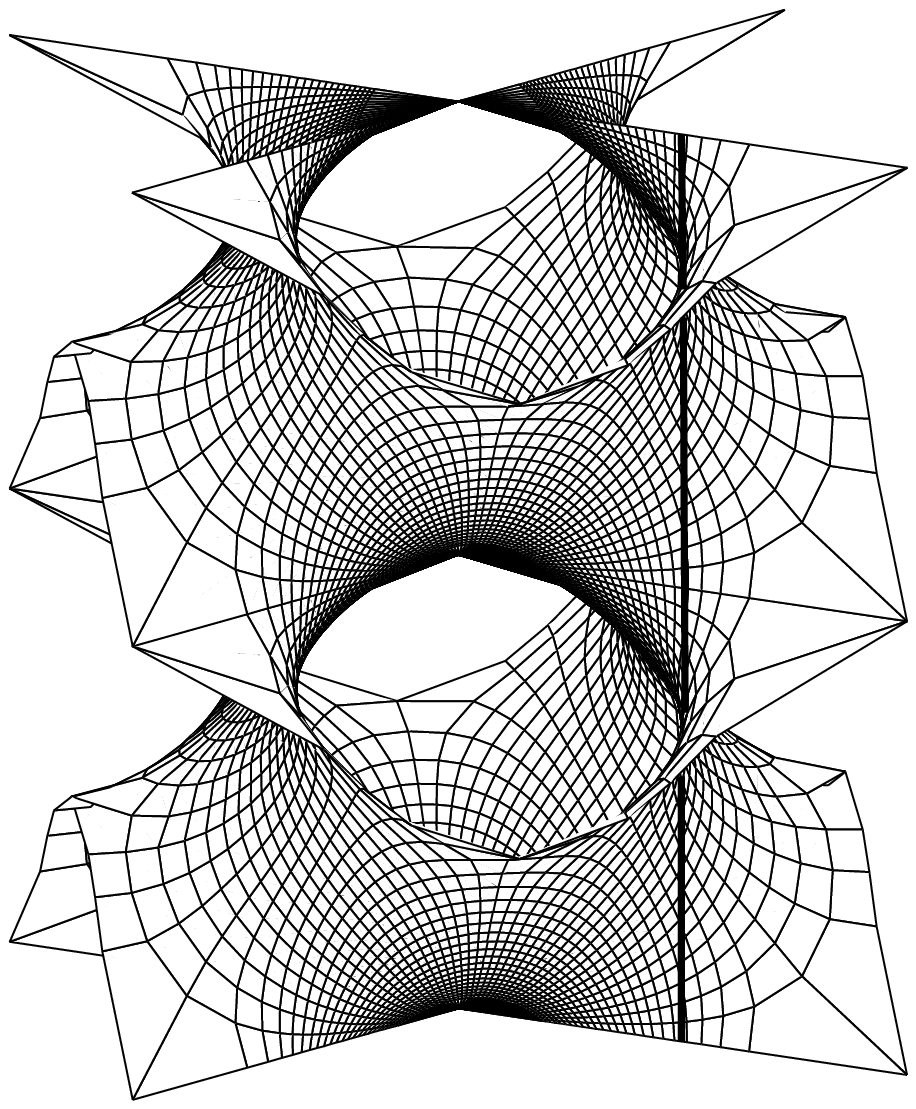}
\epsfxsize 9cm  
\epsfbox{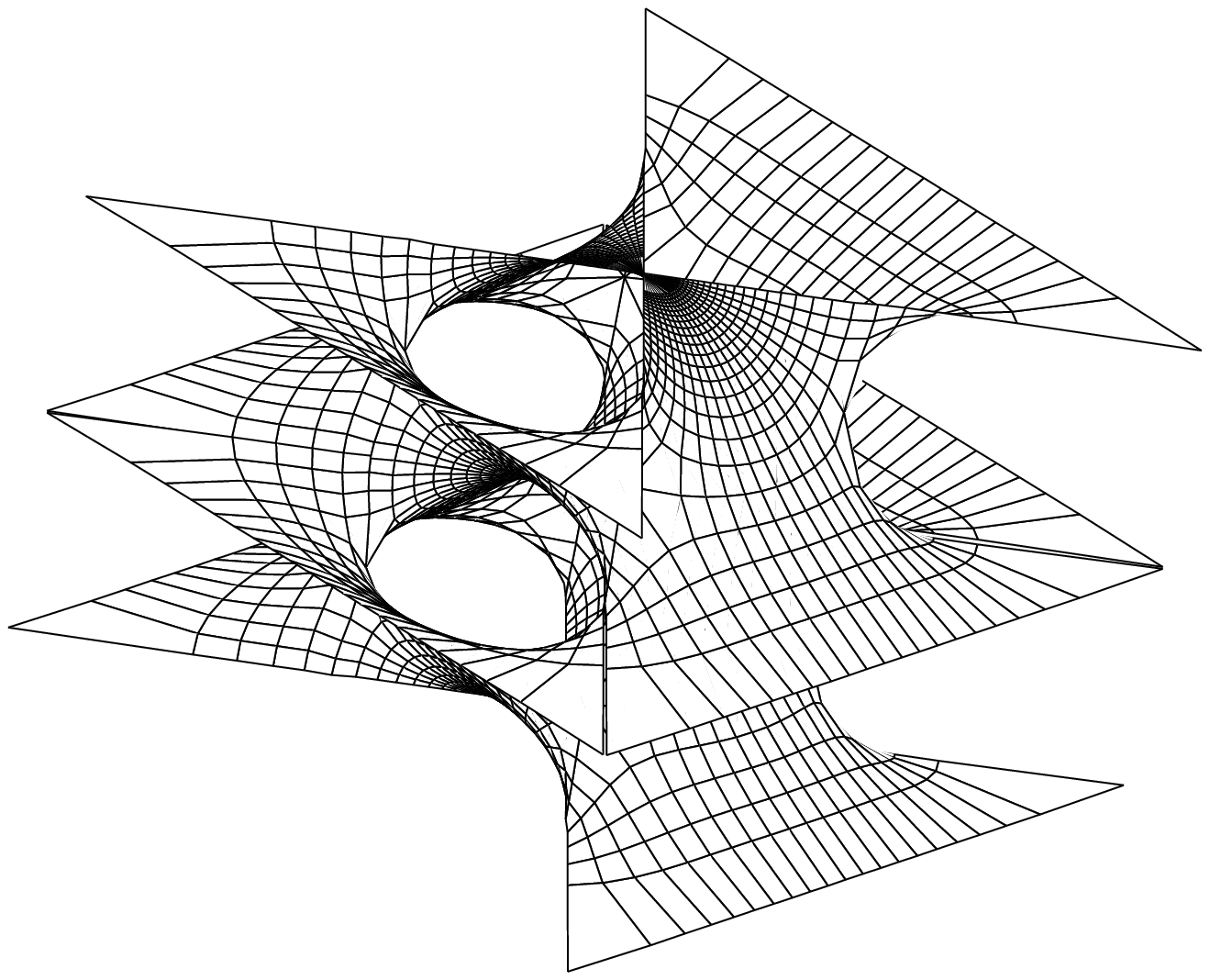}}
\centerline{(a)\hspace{7cm}(b)}  
\caption{(a) Schwarz's P-surfaces; (b) the surfaces L$_4$.} 
\end{figure} 

We are going to prove the following results:  
\\ 
\\ 
{\bf Theorem 1.1.} \it There exists a one-parameter family of triply periodic minimal surfaces in $\R^3$, of which the members are called C$_2$, and for any of them the following holds:
\\
(a) The quotient by its translation group has genus 9.
\\
(b) The whole surface is generated by a fundamental piece, which is a surface with boundary in $\R^3$. The boundary consists of four planar curves of vertical reflectional symmetry. The fundamental piece has a symmetry group generated by two vertical planes of reflectional symmetry and two line segments of 180$^\circ$-rotational symmetry.
\\
(c) By successive reflections in the boundary of the fundamental piece, and successive vertical translations, one obtains the triply periodic surface.\rm  
\\
\\
{\bf Theorem 1.2.} \it For $k=2,4$, there exists a one-parameter family of triply periodic minimal surfaces in $\R^3$, of which the members are called L$_k$, and for any of them the following holds:
\\
(a) The quotient by its translation group has genus $2k+1$.
\\
(b) The whole surface is generated by a fundamental piece, which is a surface with boundary in $\R^3$. The boundary consists of four line segments. The fundamental piece has a symmetry group generated by two vertical planes of reflectional symmetry and two line segments of 180$^\circ$-rotational symmetry. Each of these segments make an angle of $\pi/k$ with the boundary.
\\
(c) By successive rotations about the boundary of the fundamental piece, and successive vertical translations, one obtains the triply periodic surface.\rm  
\\

Sections 3 to 7 are devoted to the proof of Theorem 1.1. The proof of Theorem 1.2 follows very similar arguments and we briefly discuss it in Section 8. For this present work, the third author was supported by the grants ``Bolsa de Produtividade Cient\ih fica'' from CNPq - Conselho Nacional de Desenvolvimento Cient\ih fico e Tecnol\'ogico, and ``Bolsa de P\'os-Doutorado'' FAPESP 2000/07090-5.

\section{Preliminaries}
In this section we state some basic definitions and theorems. Throughout this work, surfaces are considered connected and regular. Details can be found in \cite{[K]}, \cite{[LpM]}, \cite{[N]} and \cite{[O]}. 
\\
\\
{\bf Theorem 2.1.} \it Let $X:R\to\E$ be a complete isometric immersion of a Riemannian surface $R$ into a three-dimensional complete flat space $\E$. If $X$ is minimal and the total Gaussian curvature $\int_R K dA$ is finite, then $R$ is biholomorphic to a compact Riemann surface $\overline{R}$ punched at a finite number of points.\rm 
\\
\\
{\bf Theorem 2.2.} (Weierstra\ss \ Representation). \it Let $R$ be a Riemann surface, $g$ and $dh$ meromorphic function and 1-differential form on $R$, such that the zeros of $dh$ coincide with the poles and zeros of $g$. Suppose that $X:R\to\E$, given by
\BE
   X(p):=Re\int^p(\phi_1,\phi_2,\phi_3),\eh\eh where\eh\eh
   (\phi_1,\phi_2,\phi_3):=\m(1/g-g,i/g+ig,2)dh,\label{W_rep}
\EE
is well-defined. Then $X$ is a conformal minimal immersion. Conversely, every conformal minimal immersion $X:R\to\E$ can be expressed as (1) for some meromorphic function $g$ and 1-form $dh$.\rm
\\
\\
{\bf Definition 2.1.} The pair $(g,dh)$ is the \it Weierstra\ss \ data \rm and $\phi_1$, $\phi_2$, $\phi_3$ are the \it Weierstra\ss \ forms \rm on $R$ of the minimal immersion $X:R\to X(R)\subset\E$.
\\
\\
{\bf Theorem 2.3.} \it Under the hypotheses of Theorems 2.1 and 2.2, the Weierstra\ss \ data $(g,dh)$ extend meromorphically on $\overline{R}$.\rm  
\\

The function $g$ is the stereographic projection of the Gau\ss \ map $N:R\to S^2$ of the minimal immersion $X$. It is a covering map of $\hat\C$ and $\int_SKdA=-4\pi$deg$(g)$. These facts will be largely used throughout this work.

\section{The surfaces $\bar{M}$ and the functions $z$} 
Consider the surface indicated in Figure 1(a). A reflection in any of its vertical planar curves of the boundary leads to a fundamental piece which represents the quotient of a triply periodic surface $M$ by its translation group. We are going to denote this quotient by $\bar{M}$. It is not difficult to conclude that it has genus 9. The fundamental domain of $\bar{M}$ is the shaded region indicated in Figure 3(a).
 
\input epsf  
\begin{figure} [ht] 
\centerline{ 
\epsfxsize 12cm 
\epsfbox{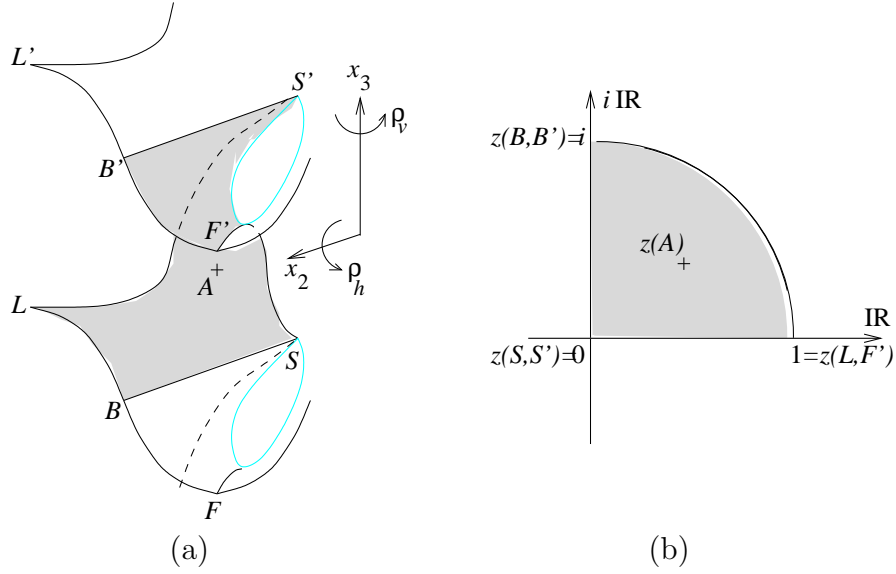}} 
\hspace{1.2in}(a)\hspace{2.3in}(b) 
\caption{(a) The fundamental domain of $M$; (b) the function $z$ on $M$.} 
\end{figure} 
 
The surface $\bar{M}$ is invariant under $180^\circ$-rotations around the directions $\vec{x}_3$ and $\vec{x}_2$. These rotations we call $r_v$ and $r_h$, respectively (see Figure 3(a)). Based on this picture, one sees that the fixed points of $r_v$ are $S,S',L,L',F,F'$ and the images of $S$ and $S'$ under the symmetries of $\bar{M}$. They sum up 8 in total. The quotients by $r_v$ and $r_h$ we call $\rho_v$ and $\rho_h$, respectively. The surface $\rho_v(\bar{M})$ is still invariant under the rotation $r_h$. In this case, the fixed points of $r_h$ will be $\rho_v(A),\rho_v(A')$ and their images under the symmetries of $\rho_v(\bar{M})$. They sum to 8 in total. Because of that,
\[
   \chi(\rho_h(\rho_v(\bar{M})))=\frac{1-9+\frac{\ds 8}{\ds 2}}{2}+\frac{8}{2}=2.
\]

Let us define $z:=\rho_h\circ\rho_v:\bar{M}\to S^2\approx\hat{\C}$, such that $z(S)=0,z(L)=1$ and $z(B)=i$. The involutions of $\bar{M}$ are induced by $\rho_v$ and $\rho_h$ on $\hat{\C}$, and since all the involutions of $\hat\C$ are given by M\"{o}bius transformations, we can conclude the following: $z(S')=0,z(F')=-z(F)=-z(L')=1$ and $z(B')=i$. By applying the symmetries of $\bar{M}$, one easily reads off the other values of $z$ at the images (under these symmetries) of $S,S',L,L',F,F',B$ and $B'$. Regarding the points $A$ and $A'$, we have $z(A)=x\in\hat\C$ such that $|x|<1$ and $Arg(x)\in(0,\pi/2)$. Consequently, $z(A')=-\bar{x}$ and one easily gets the other values of $z$ at the images of $A$ and $A'$ under the symmetries of $\bar{M}$.

\section{The $g$-function on $\bar{M}$ in terms of $z$}
First of all, observe that Jorge-Meeks' formula gives deg$(g)=9-1=8$. Let us then consider Figure 3(b). We shall have $g-g^{-1}=\infty$ if and only if $z-z^{-1}\in\{0,\infty\}$. Moreover, $g-g^{-1}=0$ if and only if $z\in\{-x,\bar{x},x^{-1},-\bar{x}^{-1},ia,ia^{-1}\}$, where $a\in(0,1)$. From this point on we introduce the following notation:
\[
   Z:=z-z^{-1},\eh X:=x^{-1}-x \eh\eh{\rm and}\eh\eh \A:=a+a^{-1}.
\] 

Based on Figure 3 it is not difficult to conclude that
\BE    
   \biggl(g-\frac{1}{g}\biggl)^2=\frac{-ic}{Z^3}\cdot(Z-i\A)^2(Z-X)(Z+\bar{X}),
\EE 
where $c$ is a positive constant. Now we define $\bar{M}$ as a member of the family of compact Riemann surfaces given by the algebraic equation (2). Later on we are going to verify that $\bar{M}$ has genus 9, indeed. But first we derive some conditions on the variables $a,x$ and $c$ in order to guarantee that $g^2=-1$ at $z=-ia^{\pm 1}$. This will be the case if
\BE
   c=\frac{\A}{\A^2+2\A Im\{X\}+|X|^2}.
\EE

Since $|X|^2=Im^2\{X\}+Re^2\{X\}$, one easily sees that $c$ is positive.

Now we analyse what happens to (2) under the map $z\to\bar{z}$. In this case we shall get $g\to i\bar{g}$ or $g\to -i\bar{g}$. Therefore 
\BE    
   \biggl(g+\frac{1}{g}\biggl)^2=\frac{-ic}{Z^3}\cdot(Z+i\A)^2(Z-\bar{X})(Z+X).
\EE 

At this point we are ready to prove that $\bar{M}$ has genus 9. The function $z$ is a four-sheet branched covering of the sphere. The values $0,\infty,\pm 1,\pm x^{\pm 1},\pm\bar{x}^{\pm 1}$ correspond to the unique branch points of $z$, all of them of order $2$, and each of these values is taken twice on $\bar{M}$. Therefore, from the Riemann-Hurwitz's formula we have
\[
   genus(\bar{M})=\frac{12\cdot(2-1)\cdot 2}{2}-4+1=9.
\] 

Now we are ready to find some relations that the parameters $a$, $c$ and $x$ will have to satisfy. These relations will make (2) and (4) consistent with the values of $g$ on the symmetry curves and lines of $\bar{M}$.

\section{Conditions on the parameters $a$, $c$ and $x$} 
Consider the curves $S'L$ and $F'S$ represented in Figure 3. The same picture shows how we have positioned our coordinate system. On the curve $S'L$ we expect that $g\in\exp(i\pi/4)\R$ and on $F'S$ one should have $g\in\exp(-i\pi/4)\R$. Let us now verify under which conditions this will really happen. 

On $S'L$ we ought to have $Re\{(g-g^{-1})^2\}=-2$. By taking $z(t)=t$, $0<t<1$, defining $T:=t-t^{-1}$ and applying it to (2) we get the following equality:  
\BE
   (g-g^{-1})^2\biggl|_{_{z(t)}}=\frac{-ic}{T^3}\cdot(T-i\A)^2(T^2-2iIm\{X\}\cdot T-|X|^2).
\EE

Therefore,
\BE
   Re\{(g-g^{-1})^2\}=-\frac{2c}{T^2}\cdot(Im\{X\}\cdot T^2+\A T^2-Im\{X\}\cdot\A^2-\A|X|^2)
\EE
on the curve $z(t)$. Since we want $Re\{(g-g^{-1})^2\}=-2$ on this curve, (6) will then give rise to the following conditions
\BE
   c=\frac{1}{\A+Im\{X\}}\eh\eh{\rm and}
\EE
\BE
   \A=-\frac{|X|^2}{Im\{X\}}.
\EE

Equation (7) can be deduced from (3) and (8) by a simple calculation. Equation (8) will restrict the definition domain of our parameters. Since $a\in(0,1)$, then $\A>2$ and by taking $x=|x|\exp(i\theta)$ one clearly sees that $Im\{X\}<0$ for $\theta\in(0,\frac{\ds\pi}{\ds 2})$. From (8) we finally get the following restriction for the $x$-variable
\BE
   Re^2\{X\}>-2Im\{X\}-Im^2\{X\}. 
\EE

Figure 4 illustrates the $X$-domain established by (9), and we recall that $|x|<1$ and $\theta\in(0,\pi/2)$. 

\input epsf  
\begin{figure} [ht] 
\centerline{ 
\epsfxsize 7cm 
\epsfbox{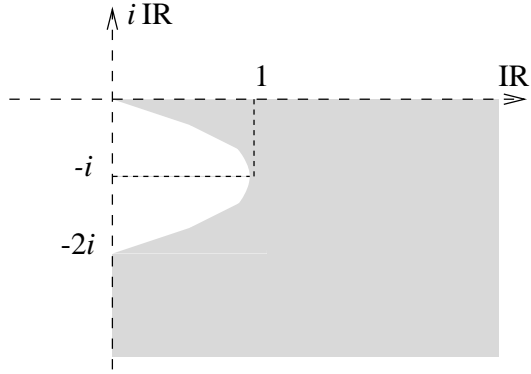}}
\caption{Definition domain of the $X$-variable satisfying (8).} 
\end{figure} 

It is not difficult to prove that (9) is equivalent to the following inequality:
\BE
   |x|<\m\biggl\{\sin\theta+\sqrt{1+3\cos^2\theta}-
       \sqrt{2\sin\theta\biggl(\sqrt{1+3\cos^2\theta}-\sin\theta\biggl)}
       \biggl\}.
\EE

Of course, the right-hand side of (10) is one of the two roots of a 2$^{\rm nd}$-degree equation. One easily proves that the other root is bigger than 1. Its inverse is exactly the right-hand side of (10), and this shows that it is positive and smaller than 1. 

Regarding our remaining restriction, namely $Re\{(g-g^{-1})^2\}=2$ on $F'S$, it is not difficult to verify that it leads to the same conditions (7) and (8). Therefore, we are now ready to write down the following table, which summarizes some special involutions of $\bar{M}$: 

\BE 
\begin{tabular}{|c|c|c|c|c|}\hline
 ${\rm symmetry}$&${\rm involution}     $&$ z{\rm-values}    $&$ g{\rm-values} $ \\ \hline\hline   
 $ SB  $&$ (g,z)\to(\bar{g},-\bar{z})   $&$ it,0<t<1         $&$ g\in\R_{+}     $ \\ \hline 
 $ BL  $&$ (g,z)\to(-\bar{g},1/\bar{z}) $&$ e^{it},\pi/2>t>0  $&$ g\in\R_{+}     $ \\ \hline 
 $ LS' $&$ (g,z)\to(i\bar{g},\bar{z})   $&$ t,1>t>0          $&$ g\in e^\q\R   $ \\ \hline   
 $ S'B'$&$ (g,z)\to(\bar{g},-\bar{z})   $&$ it,0<t<1         $&$ g\in\R_{-}     $ \\ \hline   
 $ B'F'$&$ (g,z)\to(-\bar{g},1/\bar{z}) $&$ e^{it},\pi/2>t>0  $&$ g\in\R_{-}     $ \\ \hline    
 $ F'S $&$ (g,z)\to(-i\bar{g},\bar{z})  $&$ t,1>t>0          $&$ g\in 
 e^{-\q}\R $ \\ \hline 
 $ \rho_h $&$ (g,z)\to(-1/g,z)          $&$ -x,\bar{x},x^{-1},-\bar{x}^{-1} 
 $&$ g=\pm i   $ \\ \hline      
\end{tabular} 
\EE 

The careful reader will notice that the points $(g,z)=(\pm i,-ia^{\pm 1})$ do not come out as fixed points of $\rho_h$ in (11). This is because the {\it germs} of the function $g$ at these points {\it are not} the same. This has to do with the fact that the power of $(Z+i\A)$ is a {\it multiple} of the power of $(g+g^{-1})$ in (4).

\section{The height differential $dh$ on $\bar{M}$}
Since the surface $M$ has no ends, $dh$ must be a holomorphic differential form on it. The zeros of $dh$ are exactly at the points where $g=0$ or $g=\infty$, and ord($dh)=|$ord($g)|$ at these points. They should sum up 16 in total, which is consistent with deg($dh$)=$-\chi(\bar{M})$. Let us now analyse the differential $dz$. Based on Figure 3, one sees that $dz$ has a simple zero at the points $z^{-1}(\{0,\pm 1,\pm x^{\pm 1},\pm\bar{x}^{\pm 1}\})$ and a pole of order 3 at the points $z^{-1}(\{\infty\})$. It is not difficult to conclude that
\BE
   \biggl(\frac{dh}{dz}\biggl)^2\sim\frac{(z^2-1)^2}{(z^2-x^2)(z^2-x^{-2})(z^2-\bar{x}^2)(z^2-\bar{x}^{-2})}.
\EE  

If we had a well defined square root of the function at the right-hand side of (12), then we could get an explicit formula for $dh$ in terms of $z$ and $dz$. This square root exists, indeed. By multiplying (2) and (4) it follows that
\[
   (Z-X)(Z+\bar{X})(Z-\bar{X})(Z+X)=\frac{-Z^6}{c^2(Z^2+\A^2)^2}\biggl(g^2-\frac{1}{g^2}\biggl)^2,
\]
which allows us to define
\BE
   \sqrt{(Z^2-X^2)(Z^2-\bar{X}^2)}:=\frac{iZ^3}{c(Z^2+\A^2)}\biggl(g^2-\frac{1}{g^2}\biggl).
\EE

Now we apply (13) to (12) and obtain
\BE
   dh=\frac{Z}{\sqrt{(Z^2-X^2)(Z^2-\bar{X}^2)}}\cdot\frac{dz}{z}.
\EE  

At (14) the equality sign holds because we want $Re\{dh\}=0$ on the straight line segment $SB$ (see Figure 3(a)). On this segment $z$ is pure imaginary and then we can fix both sides of (14) to be equal. Let us now verify if the symmetry curves and lines of $M$ really exist. From (11) and (14) we write down the following table 
\BE 
\begin{tabular}{|c|c|c|c|c|}\hline
 ${\rm symmetry} $&$ z{\rm-values} $&$ g\in $&$ dh(\dot{z})\in $\\ \hline\hline  $ SB  $&$ it,0<t<1        $&$ \R_{+}   $&$ i\R $\\ \hline 
 $ BL  $&$ e^{it},\pi/2>t>0 $&$ \R_{+}   $&$ \R  $\\ \hline 
 $ LS' $&$ t,1>t>0         $&$ e^\q\R   $&$ \R $\\ \hline   
 $ S'B'$&$ it,0<t<1        $&$ \R_{-}   $&$ i\R $\\ \hline   
 $ B'F'$&$ e^{it},\pi/2>t>0 $&$ \R_{-}   $&$ \R  $\\ \hline    
 $ F'S $&$ t,1>t>0         $&$ e^{-\q}\R $&$ \R  $\\ \hline      
\end{tabular} 
\EE 

From (15) it follows that $\frac{\ds dg}{\ds g}\cdot dh$ is purely imaginary on $SB$ and $S'B'$. It is real on the other paths, confirming that $M$ will have the expected symmetry curves and lines.

\section{Solution of the period problems}
The analysis of the period problems can be reduced to the analysis of the fundamental domain of our minimal immersion. If this fundamental domain is contained in a rectangular prism of $\R^3$, and if the boundary of the former is contained in the border of the latter, we shall have that the fundamental piece of our minimal surface will be free of periods. 

In order to obtain such a prism, a little reflection will show us that the following two conditions will be enough:
\\
1. The symmetry $\rho_h$ really exists in $\R^3$;
\\
2. After an orthonormal projection of the fundamental domain in the direction $x_3$, we shall have $S=S'$ and $B=B'$ (see Figure 5).

\input epsf  
\begin{figure} [ht] 
\centerline{ 
\epsfxsize 8cm 
\epsfbox{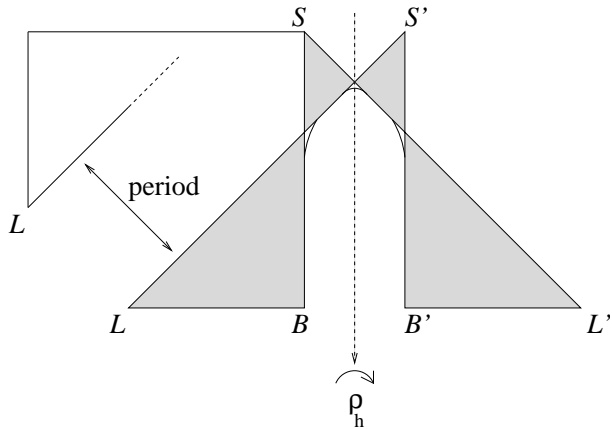}}
\caption{$x_3$-projection of the fundamental domain with an open period.} 
\end{figure} 

The first condition is easy to prove. Take a path $P\to A\to P'$ on $\bar{M}$ as indicated in Figure 6. Consider that $A\to P'$ with reversed orientation is the image of $P\to A$ under the involution $(g,z)\to(-1/g,z)$. Now we compute in $\R^3$ what happens to the coordinates of our minimal surface:  
\[
   (x_1,x_2,x_3)|_{_{(g,z)\to(-1/g,z)}}=Re\int_{P=(g_0,z_0)}^{A=(-i,x)}(\phi_1,\phi_2,\phi_3)=
\]
\[
   Re\int_{P'=(-1/g_0,z_0)}^{A=(-i,x)}(\phi_1,-\phi_2,\phi_3)=Re\int_{P=(g_0,z_0)}^{A=(-i,x)}(-\phi_1,\phi_2,-\phi_3)=(-x_1,x_2,-x_3).
\]

\input epsf  
\begin{figure} [ht] 
\centerline{ 
\epsfxsize 14cm 
\epsfbox{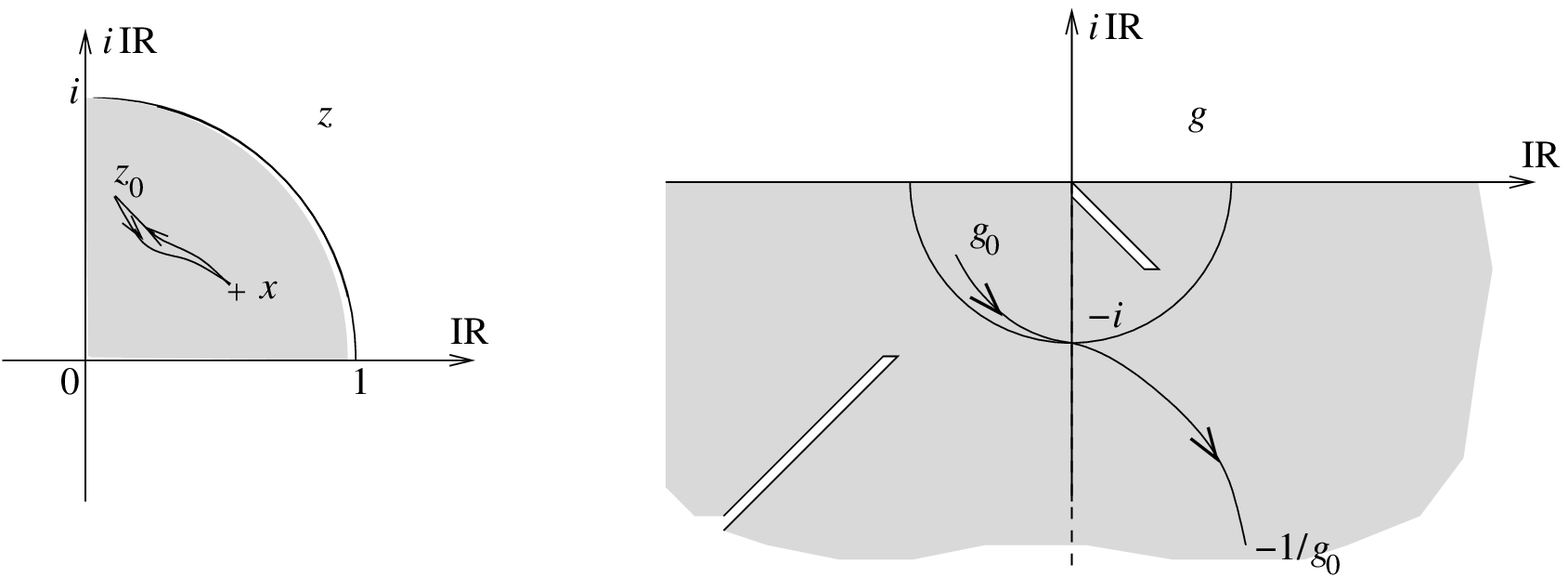}}
\caption{The path $P\to A\to P'$ on $\bar{M}$.} 
\end{figure} 

Therefore, our minimal surface is {\it really} invariant under 180$^\circ$-rotations around the $x_2$-axis. This proves the existence of the symmetry $\rho_h$ of our initial assumptions.

Now we are ready to deal with the second condition. Consider Figure 5 with the segments $SB$ and $BL$ on it. The period will be zero if and only if these segments have the same length, or equivalently
\BE
   Re\int_{SB}\phi_2=Re\int_{BL}\phi_1.
\EE    

On $SB$ we can take $Z(t)=it$, $2<t<\infty$. This implies that $\frac{\ds dz}{\ds z}=-\frac{\ds dt}{\ds\sqrt{t^2-4}}$. From (4) and (14)
we have
\BE
   \phi_2\biggl|_{_{Z(t)=it}}=\frac{c^\m(t+\A)}{t^\m(t^2-2Im\{X\}\cdot t+|X|^2)^\m}\cdot\frac{dt}{\sqrt{t^2-4}}.
\EE

On $BL$ we can take $z(t)=\exp(it)$, $0<t<\frac{\ds\pi}{\ds 2}$. From (2) and (14) it follows that
\BE
   \phi_1\biggl|_{_{z(t)=\exp(it)}}=\frac{1}{\sqrt{2}}\cdot\frac{c^\m(\A-2\sin t)}{(4\sin^2t+4Im\{X\}\cdot\sin t+|X|^2)^\m}\cdot\frac{dt}{\sqrt{\sin t}}.
\EE

Now define $I_1:=\frac{1}{\sqrt{2c}}\int_{BL}\phi_1$ and $I_2:=\frac{1}{\sqrt{2c}}\int_{SB}\phi_2$. For $I_1$ apply the change of variables $u^2=\sin t$ and for $I_2$, $t=2u^{-2}$. A simple reckoning will lead to the following equalities:
\BE
   I_1=\int_{0}^{1}\frac{\A-2u^2}{(4u^4+4Im\{X\}\cdot u^2+|X|^2)^\m}\cdot\frac{du}{\sqrt{1-u^4}}
\EE
and
\BE
   I_2=\int_{0}^{1}\frac{2+\A u^2}{(4-4Im\{X\}\cdot u^2+|X|^2u^4)^\m}\cdot\frac{du}{\sqrt{1-u^4}}.
\EE 

The next proposition will solve the period problem given by (16):
\\
\\
{\bf Proposition 1.} \it For any fixed positive value of $Re\{X\}$ one has that the following limit exists and is positive
\[
   \lim_{Im\{X\}\to 0}(-Im\{X\})\cdot(I_1-I_2).
\]
For $Im\{X\}=-1$ we have that $\Lim{\A\to 2}{(I_1-I_2)}$ exists and is negative.\rm
\\
\\
{\bf Proof.} By recalling (8), a simple reckoning will show that
\BE
   \lim_{Im\{X\}\to 0}(-Im\{X\})\cdot I_1=\int_{0}^{1}\frac{Re^2\{X\}}{(4u^4+Re^2\{X\})^\m}\cdot\frac{du}{\sqrt{1-u^4}}
\EE
and
\BE
   \lim_{Im\{X\}\to 0}(-Im\{X\})\cdot I_2=\int_{0}^{1}\frac{Re^2\{X\}\cdot u^2}{(4+Re^2\{X\}\cdot u^4)^\m}\cdot\frac{du}{\sqrt{1-u^4}}.
\EE

Since
\[
   \frac{u^2}{(4+Re^2\{X\}\cdot u^4)^\m}<\frac{1}{(4u^4+Re^2\{X\})^\m}
\]
for every $Re\{X\}>0$ and $u\in(0,1)$, from (21) and (22) it follows the first
assertion of Proposition 1. 

By fixing $Im\{X\}=-1$ and recalling (9), the convergence $\A\to 2$ is equivalent to $Re\{X\}\to 1$. This means that $X$ approaches the point $1-i$ indicated in Figure 4. An easy calculation will give us
\BE
   \lim_{\A\to 2} I_1=\sqrt{2}\int_{0}^{1}(2u^4-2u^2+1)^{-\m}\cdot\biggl(\frac{1-u^2}{1+u^2}\biggl)^\m du
\EE
and
\BE
   \lim_{\A\to 2} I_2=\sqrt{2}\int_{0}^{1}(2+2u^2+u^4)^{-\m}\cdot\biggl(\frac{1+u^2}{1-u^2}\biggl)^\m du.
\EE

The integrand of (23) can be rewritten as
\BE
   (u^4+(1-u^2)^2)^{-\m}\cdot\biggl(\frac{1-u^2}{1+u^2}\biggl)^\m=\biggl[\frac{u^4}{(1-u^2)^2}+1\biggl]^{-\m}\cdot\frac{1}{\sqrt{1-u^4}},
\EE
while one rewrites the integrand of (24) as 
\BE
   (u^4+(1+u^2)^2)^{-\m}\cdot\biggl(\frac{1+u^2}{1-u^2}\biggl)^\m=\biggl[\frac{u^4}{(1+u^2)^2}+1\biggl]^{-\m}\cdot\frac{1}{\sqrt{1-u^4}}.
\EE

Since 
\[
   \frac{u^4}{(1-u^2)^2}>\frac{u^4}{(1+u^2)^2}
\]
for every $u\in(0,1)$, it follows the last assertion of Proposition 1.\hfill $\square$
\ \\

Proposition 1 provides a family of triply periodic surfaces of which a member is exemplified in Figure 1(a). By looking at Figure 4, this family can be represented by the values of $X$ which belong to a curve $\CC$ contained in the shaded region. All members of this family will have {\it only} three periods, as suggested by Figure 1(a). Nevertheless, a priori there might be some non-embedded members, but it will not be the case. This is the subject of our next section.

\section{Embeddedness of the triply periodic surfaces}
From now on we shall denote our triply periodic surfaces by $M_X$, where $X\in\CC$. Figure 6 shows that the projection of the unitary normal on a fundamental domain of $M_X$ is contained in the lower hemisphere of $\hat\C$. This means that $(x_1,x_3)$ is an immersion of $\D:=\{z\in\C:|z|<1\eh{\rm and}\eh 0<Arg(z)<\pi/2\}$ in $\R^2$. The next picture shows a possible image of this map in $\R^2$:   

\input epsf  
\begin{figure} [ht] 
\centerline{ 
\epsfxsize 10cm  
\epsfbox{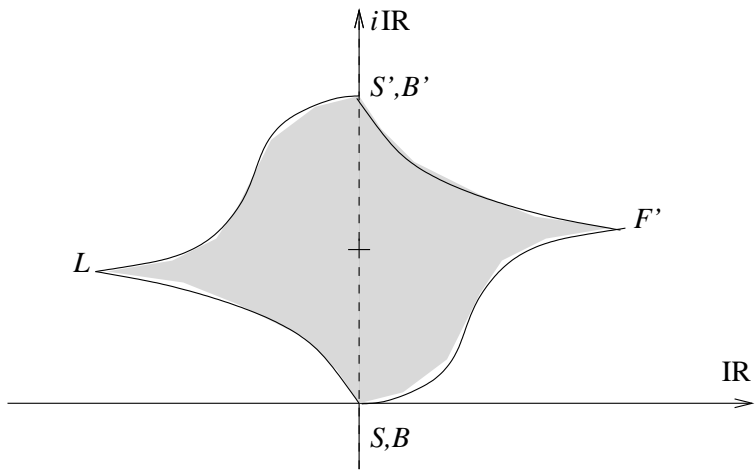}}  
\caption{A possible $x_2$-projection of the fundamental domain on $x_1Ox_3$.} 
\end{figure} 

It is not difficult to prove that the contour of the shaded region in Figure 7 is a monotone curve. The $x_1$-coordinate of the curve $BL$ is given by the integral of $-\phi_1$ as in (18). The integrand is clearly positive, hence this stretch is monotone. Regarding $LS'$, where we can take $Z(t)=t$, $0>t>-\infty$, a simple reckoning gives us
\[
   dh\biggl|_{_{Z(t)=t}}=\frac{t}{\sqrt{t^4-Re\{X^2\}t^2+|X|^4}}\cdot\frac{dt}{\sqrt{t^2+4}}.
\]

Hence, the stretch $LS'$ is also monotone. By using the symmetry $\rho_h$, it follows that the whole contour indicated in Figure 7 is a monotone curve. Since the 3$^{\rm rd}$ coordinate of $BL$ is increasing, the projections $BL$ and $LS'$ will intersect only at the point $L$. Nevertheless, it can happen that the projection $LS'$ crosses $B'F'$. If we prove that this is not the case, the contour will have no self-intersections. The shaded region will then be simply connected, and we shall conclude that the fundamental domain is a graph, hence embedded. 

But even so, it can happen that the expanded triply periodic surface will not be embedded. We do not know whether the curve $LS'$ crosses the $x_3$-axis or not. A little reflection will show that, if $g$ does not take the value $-\exp(i\pi/4)$ along $LS'$, then this curve does not intersect the vertical axis. Consequently, its projection will not intersect $B'F'$. In this case, since the triply periodic surface is expanded horizontally by reflections only, and vertically by rotations only, the whole surface will then be embedded.

By using the maximum principle, if we find an embedded member of our family in
the curve $\CC$, the whole family will then consist of embedded surfaces. The
following proposition gives us such a member and will conclude this section:
\\
\\
{\bf Proposition 2.} \it There is an $X\in -i+(1,2\sqrt{2})$ such that $X\in\CC$
and $M_X$ is embedded.\rm
\\
\\
{\bf Proof.} We shall prove that $g\ne-\exp(i\pi/4)$ along $LS'$, for any $X\in-i+(1,2\sqrt{2})$. Moreover, $(I_1-I_2)|_{_{X=2\sqrt{2}-i}}$ will be positive. These two facts together with Proposition 1 will conclude Proposition 2. 
\\

By recalling (5), we would have $g=-\exp(i\pi/4)$ for some $T\in(-\infty,0)$ if and only if
\BE
   (T^2-\A^2)(T^2-|X|^2)=4\A Im\{X\}\cdot T^2.
\EE

Equation (27) will not be fulfilled by any $T^2\in(0,\infty)$ providing 
\[
   |\A^2+|X|^2+4\A Im\{X\}|<2\A |X|,
\]
or equivalently 
\[
   -\frac{Re\{X\}}{Im\{X\}}<2\sqrt{2}.
\]   

We have fixed $Im\{X\}=-1$, hence $g\ne-\exp(i\pi/4)$ along $LS'$ for any $X\in-i+(1,2\sqrt{2})$. Let us now verify that $(I_1-I_2)|_{_{X=2\sqrt{2}-i}}>0$. From (19) we have 
\BE
   I_1|_{_{X=2\sqrt{2}-i}}=
   \int_{0}^{1}\frac{9-2u^2}{(4u^4-4u^2+9)^\m}\cdot\frac{du}{\sqrt{1-u^4}}
\EE
and from (20) it follows that 
\BE
   I_2|_{_{X=2\sqrt{2}-i}}=\int_{0}^{1}\frac{2+9u^2}{(4+4u^2+9u^4)^\m}\cdot\frac{du}{\sqrt{1-u^4}}.
\EE 

But 
\BE
   \frac{9-2u^2}{(4u^4-4u^2+9)^\m}>3-\frac{2}{3}u^2,\eh\forall\eh u\in(0,1),
\EE
and if we define $a:=1-11/\sqrt{17}$ it is possible to prove that 
\BE
   \frac{2+9u^2}{(4+4u^2+9u^4)^\m}<au^2-2au+1,\eh\forall\eh u\in(0,1).
\EE

But 
\BE
   \tilde{I}_1:=\int_{0}^{1}\frac{(3-2u^2/3)du}{\sqrt{1-u^4}}=\frac{3}{4}B(\frac{1}{4},\m)-\frac{1}{6}B(\frac{3}{4},\m),
\EE
and 
\BE
   \tilde{I}_2:=\int_{0}^{1}\frac{(au^2-2au+1)du}{\sqrt{1-u^4}}=\frac{a}{4}B(\frac{3}{4},\m)-\frac{a}{2}B(\m,\m)+\frac{1}{4}B(\frac{1}{4},\m).
\EE

Now we use $B(m,n)=\Gamma(m)\Gamma(n)/\Gamma(m+n)$, $\Gamma(\frac{1}{4})=3,625600...$, $\Gamma(\m)=\sqrt{\pi}$ and $\Gamma(\frac{3}{4})=1,225417...$ in order to conclude that
\BE
   \tilde{I}_1>\tilde{I}_2.
\EE 

Together with (28-33), (34) shows that $(I_1-I_2)|_{_{X=2\sqrt{2}-i}}$ is positive.\hfill $\square$

\section{The surfaces L$_{2,4}$}
In order to prove Theorem 1.2, one follows very similar ideas already explained in Sections 3 to 7. For the surfaces L$_2$, consider Figures 8(a) and 8(b). The fundamental piece $\bar{M}$ has genus 5, and $Ox_2$ passes through point $A$. The piece is invariant under $r_v$ and $r_h$, with quotient functions $\rho_v$ and $\rho_h$, respectively. 

\input epsf  
\begin{figure} [ht] 
\centerline{ 
\epsfxsize 15cm 
\epsfbox{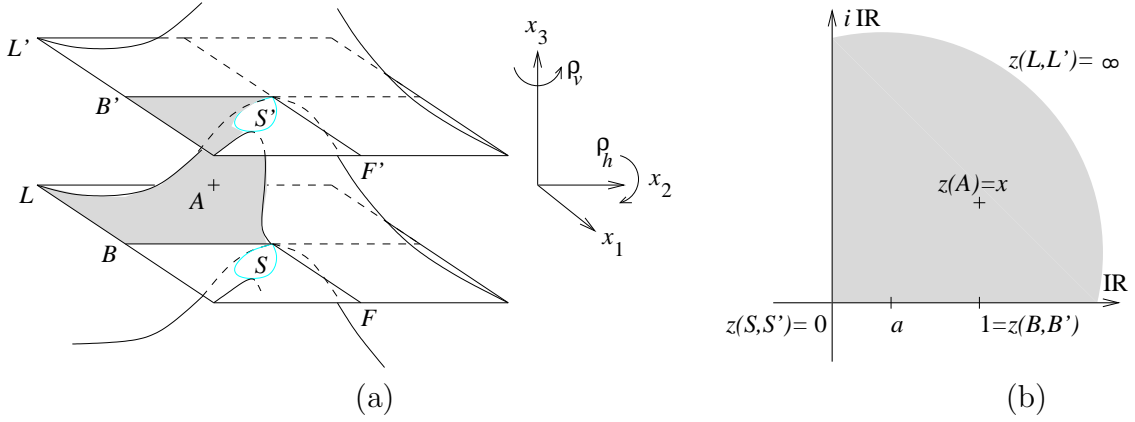}} 
\hspace{1.6in}(a)\hspace{3.2in}(b) 
\caption{(a) The fundamental domain of $M$; (b) the function $z$ on $M$.} 
\end{figure} 

Since 
\[
   \chi(\rho_h(\rho_v(\bar{M})))=\frac{1-5+\frac{\ds 8}{\ds 2}}{2}+\frac{4}{2}=2,
\]
we may define $z:=\rho_h\circ\rho_v:\bar{M}\to S^2\approx\hat{\C}$, such that $z(S)=0,z(B)=1$ and $z(L)=\infty$. The symmetries imply $z(S')=0,z(B')=1$ and $z(L')=\infty$, whereas $z(A)$ is a certain complex $x$ in the first quadrant. Moreover, there is a point in the segment $BS$ at which $g=1$. After analysing the divisors of $z$ and $g$ on $\bar{M}$, together with the behaviour of the unitary normal on symmetry curves and lines, we get  
\BE
   \biggl(g+\frac{1}{g}\biggl)^2=\frac{1/a-a}{|x-a|^2}\cdot\frac{(z-x)(z-\bar{x})(z+a)^2}{z(1-z^2)}.
\EE
Since there is a point in the segment $FS$ at which $g=-i$, we should also have
\BE
   \biggl(g-\frac{1}{g}\biggl)^2=\frac{1/a-a}{|x-a|^2}\cdot\frac{(z+x)(z+\bar{x})(z-a)^2}{z(1-z^2)}.
\EE
In order to have equivalence between (35) and (36), a necessary and sufficient condition is $\A=a+a^{-1}=\frac{|x|^2+1}{Re\{x\}}$. Now, it is easy to get
\BE
   dh=\frac{idz}{\sqrt{(z^2-x^2)(z^2-\bar{x}^2)}},
\EE 
with a well-defined square root in the denominator. One checks the assumed symmetries the same way we did in (11) and (15). The unique period problem is again (16), which can be visualised again by Figure 5. Therefore, (16) is equivalent to $J_1=J_2$, where
\BE
   J_1=\int_0^1\frac{(t+a)dt}{|t+x|\sqrt{t(1-t^2)}}   
\EE
and
\BE
   J_2=\int_1^\infty\frac{(t-a)dt}{|t-x|\sqrt{t(t^2-1)}}.   
\EE
The change $t\mapsto1/t$ for $J_2$ makes clear that $J_1<J_2$ ($J_1>J_2$) providing $R_1<R_2$ ($R_1>R_2$), where $R_1=\frac{t+a}{1-at}$ and $R_2=\biggl|\frac{t+x}{1-xt}\biggl|$, $0<t<1$. On the one hand, for a fixed $r=Re\{x\}>1$, if $Im\{x\}\to 0$ then $a\to1/r$, and consequently $R_1<R_2$. On the other hand, by fixing $Im\{x\}$ and letting $Re\{x\}\to 0$, then $a\to 0$ and so $R_1>R_2$. In this case, notice that the singularity at $t=1$ of both integrands in (38) and (39) is easily removable with a change of variables. This means, no matter we have $R_1|_{t=1}=R_2|_{t=1}$, it still holds $J_1>J_2$. 

For the surfaces L$_4$, consider Figures 9(a) and 9(b). The fundamental piece $\bar{M}$ has genus 9, and $Ox_2$ passes through point $A$. The piece is invariant under $r_v$ and $r_h$, with quotient functions $\rho_v$ and $\rho_h$, respectively. We shall have $g-g^{-1}=\infty$ if and only if $z\in\{\pm i,0,\infty\}$. Moreover, $g-g^{-1}=0$ if and only if $z\in\{-x,\bar{x},-x^{-1},\bar{x}^{-1},ia,-ia^{-1}\}$, where $a\in(0,1)$. 

\input epsf  
\begin{figure} [ht] 
\centerline{ 
\epsfxsize 15cm 
\epsfbox{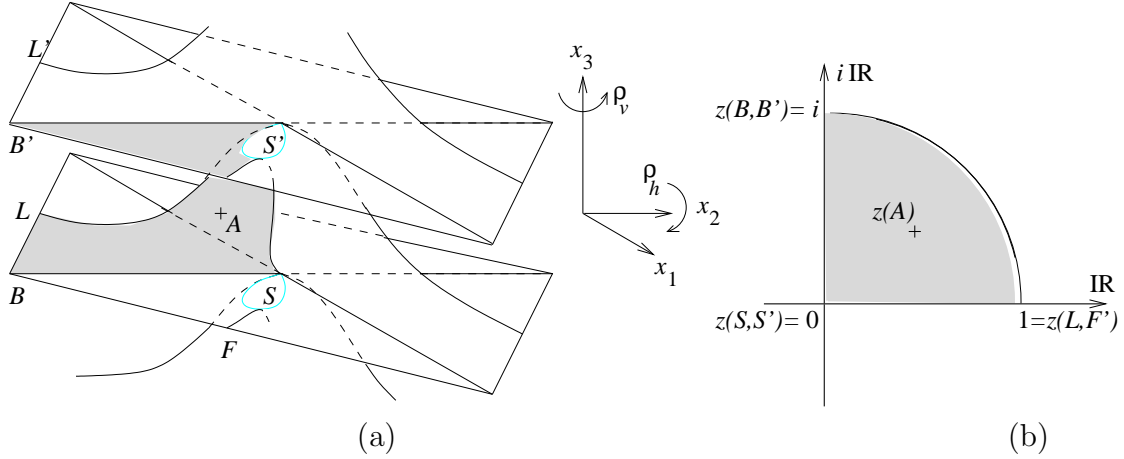}} 
\hspace{1.6in}(a)\hspace{3.2in}(b) 
\caption{(a) The fundamental domain of $M$; (b) the function $z$ on $M$.} 
\end{figure} 

From this point on we re-define:
\[
   Z:=z^{-1}+z,\eh X:=x^{-1}+x \eh\eh{\rm and}\eh\eh \A:=a^{-1}-a.
\] 

Based on Figure 9 it is not difficult to conclude that
\BE    
   \biggl(g-\frac{1}{g}\biggl)^2=\frac{ic}{Z^3}\cdot(Z+i\A)^2(Z+X)(Z-\bar{X}),
\EE 
where $c$ is given by (3) again. Moreover, (40) is equivalent to
\BE    
   \biggl(g+\frac{1}{g}\biggl)^2=\frac{ic}{Z^3}\cdot(Z-i\A)^2(Z+\bar{X})(Z-X),
\EE 
Similar arguments as in Section 5 will give again (7) and (8), but unlike (9) there is {\it no} restriction now. Regarding $dh$, it still holds (14), but unlike Figure 5 the period problem is now illustrated by Figure 10.    

Integrals $I_1$ and $I_2$ are again given by (19) and (20), but now the period is solved when $2I_1=I_2$. This will come with
\\
\\
{\bf Proposition 3.} \it For any fixed positive value of $Re\{X\}$ one has that the following limit exists and is positive
\[
   \lim_{Im\{X\}\to 0}(-Im\{X\})\cdot(2I_1-I_2).
\]
For $Im\{X\}=-1$ we have that $\Lim{\A\to 0}{(I_1-2I_2)}$ exists and is negative.\rm

The proof of Proposition 3 is quite similar to the proof of Proposition 1, and so we shall omit it here. The arguments for the embeddedness of $L_{2,4}$ are even easier than the ones used in Section 7 for $C_2$, because now the contours are given by {\it four} straight line segments and {\it two} curves, pairwise congruent.
\eject
\input epsf  
\begin{figure} [ht] 
\centerline{ 
\epsfxsize 8cm 
\epsfbox{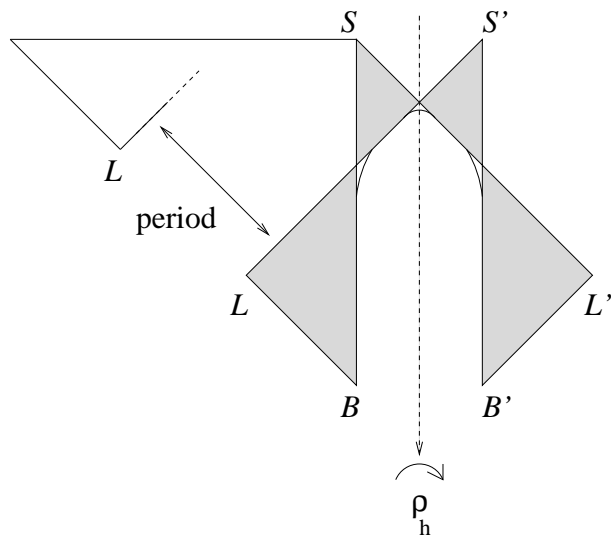}}
\caption{$x_3$-projection of the fundamental domain with an open period.} 
\end{figure} 

\end{document}